\documentclass[12pt]{article} \textwidth=6in \oddsidemargin=0in
\usepackage{amssymb} 
\begin{document}\def\ov{\over} 
\newcommand{\C}[1]{{\cal C}_{#1}} \def\inv{^{-1}}\def\be{\begin{equation}} \def\ep{\varepsilon} \newcommand{\xs}[1]{\x_{\s(#1)}}
\def\ee{\end{equation}}\def\x{\xi}\def\({\left(} \def\){\right)} \def\iy{\infty}  \def\ld{\ldots} \def\Z{\mathbb Z} \def\cd{\cdots}
 \def\P{\mathbb P}\def\bks{\backslash} \def\t{\tau} \def\s{\sigma}
 \def\r{\rho} \def\l{\ell} \def\ph{\varphi} \def\dl{\delta} \def\ps{\psi} \def\e{\eta} \def\phy{\ph_\iy} \def\z{\zeta} \def\m{\mu} 
\def\g{\gamma} \def\bs{\bar{s}} \newcommand{\br}[2]{\left[{#1\atop #2}\right]} \def\sp{\vspace{2ex}} \def\la{\lambda}
\def\noi{\noindent} \def\T{\mathcal{T}} \def\G{\Gamma} \def\a{\alpha}
\def\KA{K_{\rm Airy}} \def\A{{\rm Ai}} 

\hfill  October 24, 2009

\begin{center}{\bf \large On ASEP with Step Bernoulli Initial Condition}\end{center}

\begin{center}{\large\bf Craig A.~Tracy}\\
{\it Department of Mathematics \\
University of California\\
Davis, CA 95616, USA\\
email: tracy@math.ucdavis.edu}\end{center}

\begin{center}{\large \bf Harold Widom}\\
{\it Department of Mathematics\\
University of California\\
Santa Cruz, CA 95064, USA\\
email: widom@ucsc.edu}\end{center}

\begin{abstract}
This paper extends earlier work on ASEP to the case of step Bernoulli initial condition. The main results are a representation in terms of a Fredholm determinant for the probability distribution of a fixed particle, and asymptotic results which in particular establish KPZ universality for this probability in one regime. (And, as a corollary, for the current fluctuations.) 
\end{abstract}

\begin{center}{\bf I. Introduction and Statement of Results}\end{center}

This paper extends earlier work on ASEP \cite{TW1,TW2,TW3,TW4} to the case of step Bernoulli initial condition: each site in $\Z^+$, independently of the others, is initially occupied with probability $\r,\ 0<\r\le 1$; all other sites are initially unoccupied. The main results are a representation in terms of a Fredholm determinant for the probability distribution of a fixed particle, and asymptotic results which in particular establish \textit{KPZ universality} (see, e.g., \cite{Sp}) for this probability in one regime. (And, as a corollary, for the current fluctuations.) 

To state the results we introduce some notation, retaining as much as possible the notation of the cited papers. We denote by $x_m(t)$ the position of the $m$th particle from the left at time $t$, and by $\P(x_m(t)\le x)$ its probability distribution. We assume throughout that $q\ne0$ and use the notation $\t=p/q$.

The operator appearing in the Fredholm determinant representation has kernel
\be K(\x,\,\x')=q\,{{\x}^x\,e^{\ep(\x)t}\ov p+q\x\x'-\x}\,{\r\,(\x-\t)\ov\x-1+\r\,(1-\t)},\label{K}\ee
where
\[\ep(\x)=p\,\x\inv+q\,\x-1.\]
(Notice that when $\r=1$, ASEP with deterministic initial condition, the last factor in the kernel does not appear.) It acts on functions on $\C{R}$, a circle with center zero and radius $R$, by
\[f(\x)\to\int_{\C{R}}K(\x,\,\x')\,f(\x')\,d\x',\ \ \ (\x\in \C{R}).\footnote{All contour integrals are to be given a factor $1/2\pi i$, and are described counterclockwise unless stated otherwise.}\]
The radius $R$ is assumed so large that the denominators $p+q\x\x'-\x$ and $\x-1+\r\,(1-\t)$ are nonzero on and outside the contour.
\sp

\noi{\bf Theorem 1}. For $p,\, q\ne0$ we have
\be\P\left(x_m(t)\le x\right)=\int {\det(I-\la K)\ov \prod_{k=0}^{m-1}(1-\la\,\t^k)}\, {d\la\ov \la},\label{P}\ee
where  the integral is taken over a contour enclosing the singularities of the integrand at $\la=0$ and $\la=\t^{-k}\ (k=0,\ldots,m-1)$ . 
 
For the asymptotics\footnote{Although \cite{TW3} contained asymptotics in three regimes as $t\to\iy$, with $x$ and $m$ fixed, with $m$ fixed and $x\to\iy$, and with $m,\,x\to\iy$, we only state the extensions of the last one.} we assume $0\le p<q$, so there is a drift to the left and TASEP is a special case. Three distribution functions will arise. The first is the Gaussian
\[G(s)={1\ov\sqrt{2\pi}}\int_{-\iy}^s e^{-z^2/2}\,dz.\]
For the others we define as usual
\[\KA(x,\,y)=\int_0^\iy \A(x+z)\,\A(y+z)\,dz,\]
\[F_2(s)=\det\,(I-\KA)\ \ \ \rm {on}\ (s,\,\iy).\]
The distribution function $F_1(s)^2$ (the $^2$ is not a footnote) is given by the analogous determinant where $\KA(x,\,y)$ is replaced by
\[\KA(x,\,y)+\A(x)\,\int_{-\iy}^y\A(z)\,dz.\footnote{That $F_1(s)^2$ is the square of a distribution function from random matrix theory (see, e.g., \cite{FS}) is irrelevant. It arises here as the Fredholm determinant of the rank one perturbation of the Airy kernel.}\]
    
We use the notations $\g=q-p$ and, for $\s>0$,
\[ \s=m/t,\ \ \ c_1=-1+2\sqrt{\s},\ \ \ c_2=\s^{-1/6}\,(1-\sqrt{\s})^{2/3},\]
\[c_1'=\r\inv\s+\r-1,\ \ \ c_2'=\r\inv((1-\r)\,(\s-\r^2))^{1/2}.\]

\noi{\bf Theorem 2}. When $0\le p<q$,
\[\lim_{t\to\iy}\P\({x_m(t/\g)-c_1\,t\ov c_2\,t^{1/3}}\le s\)=F_2(s)\ \ \ {\rm when}\ 0<\s<\r^2,\]
\[\lim_{t\to\iy}\P\({x_m(t/\g)-c_1\,t\ov c_2\,t^{1/3}}\le s\)=F_1(s)^2\ \ \ {\rm when}\ \s=\r^2,\ \ \r<1,\]
\[\lim_{t\to\iy}\P\({x_m(t/\g)-c_1'\,t\ov c_2'\,t^{1/2}}\le s\)=G(s)\ \ \ {\rm when}\ \s>\r^2,\ \ \r<1.\]
These hold uniformly for $s$ in a bounded set, and parts one and three hold uniformly for $\s$ in a compact subset of its domain. Part two holds more generally for $\s=\r^2+o(t^{-1/3})$.
\sp

The \textit{total current} $\T$ at position $x$ at time $t$ is defined by
\[\T(x,t):=\textrm{ number of particles}\>\le x\>\>\textrm{at time}\>\> t.\]
We use the notations
\[v=x/t,\ \ \ a_1=(1+v)^2/4,\ \ \ a_2=2^{-4/3}(1-v^2)^{2/3},\]
\[a_1'=\r v+\r(1-\r),\ \ \  a_2'=(\r(1-\r)(v+1-2\r))^{1/2}.\]

\noi{\bf Theorem 3.} When $0\le p<q$,
\[\lim_{t\to\iy}\P\( {\T(x,t/\g)-a_1\,t \ov a_2\,t^{1/3} }\le s\)=1-F_2(-s)\ \ \ {\rm when}\ -1<v<2\r-1,\]
\[\lim_{t\to\iy}\P\( {\T(x,t/\g)-a_1\,t \ov a_2\,t^{1/3} }\le s\)=1-F_1(-s)^2\ \ \ {\rm when}\ v=2\r-1,\ \ \r<1,\]
\[\lim_{t\to\iy}\P\( {\T(x,t/\g)-a_1'\,t \ov a_2'\,t^{1/2} }\le s\)=1-G(-s)\ \ \ {\rm when}\ v>2\r-1,\ \ \r<1.\]
These hold uniformly for $s$ in a bounded set, and parts one and three hold uniformly for $v$ in a compact subset of its domain. Part two holds more generally for $v=2\r-1+o(t^{-1/3})$.
\sp

For the results in \cite{TW1,TW2} for step initial condition we began with a formula for the probability distribution with initial configuration $[1,\ld,N]$ and passed to the limit as $N\to\iy$. Here, for step Bernoulli initial condition, we first average over all initial configurations in this interval. A crucial ingredient is a new  combinatorial identity (see \S III), generalizing one in \cite{TW1}, which allows us to turn an infinite series representation for the probability distribution into the Fredholm determinant representation.

For the case of TASEP the limit formulas we obtain here (and others) were conjectured in \cite{PS} and proved in \cite{BC}.
We recommend to the reader the Introduction of \cite{BC} where a clear and detailed explanation is given of the bigger picture relating asymmetric exclusion processes to KPZ Universality and
the stochastic Burgers equation.  Ideally one would like to prove
limit theorems at this continuum level; but see \cite{BQS} (and references therein) for the difficulties of this approach.

The transition from $F_2$ to Gaussian passing through $F_1^2$ is not a new phenomenon.  It has occured also in the Corner Growth model \cite{BR}, the PNG model \cite{IS}, random matrices \cite{BBP} and in TASEP \cite{PS, BFS, BC}. In these earlier instances
the transition was more or less indirectly related to an underlying determinantal structure.  In our analysis this transition can be traced back to the pole introduced into the kernel  (\ref{K}) by the factor
$\r (\x-\t)/(\x-1+\r(1-\t))$  (which is not present for step initial condition) and its position relative to the contours through the saddle point. Since ASEP is not a determinantal process, this suggests that the $F_2\rightarrow F_1^2\rightarrow G$ transition is a universal phenomenon.

\begin{center}{\bf II. Preliminary Formula for \boldmath{$\P(x_m(t)\le x)$}}\end{center}

We begin with Theorem 5.2 of \cite{TW1} which gives a formula for $\P_Y(x_m(t)\le x)$, the probability in ASEP when the initial configuration is a (deterministic) finite set $Y$. To state it we need even more notation, which here will be slightly different than in \cite{TW1}. Define 
\[I(k,\x)=I(k,\x_1,\ld,\x_k)=\prod_{i<j}{\x_j-\x_i\ov p+q\x_i\x_j-\x_i}\;(1-\prod_i\x_i)\,\prod_i{\x_i^{x-1}\,e^{\ep(\x_i)t}\ov 1-\x_i},\]
where all indices in the products run from 1 to $k$,
and then for a set $S=\{s_1,\ld,s_k\}$ with $s_1<\cd<s_k$ define
\[I(S,\x)=I(k,\x)\,\;\prod_{i}\x_i^{-s_i}.\]
Notice that both $I(k,\x)$ and $I(S,\x)$ contain both $x$ and $t$, although they are not displayed in the notation.

The $\t$-binomial coefficient is defined by
\[\br{N}{n}_\t={(1-\t^N)\,(1-\t^{N-1})\cdots (1-\t^{N-n+1})\ov (1-\t)\,(1-\t^2)\cdots (1-\t^n)}.\]

Finally, given two sets $U$ and $V$ we define
\[\s(U,\,V)=\#\{(u,\,v): u\in U,\ v\in V,\ {\rm and}\ u\ge v.\}.\]

Theorem 5.2 of \cite{TW1} is the formula
\be\P_Y(x_m(t)=x)=\sum_{k\ge1}c_{m,k}\,\sum_{{S\subset Y\atop |S|=k}}\,\t^{\s(S,\,Y)}\,\int_{|\x_i|=R} I(x,S,\x)\,d^{|S|}\x,\label{PY}\ee
where
\be c_{m,k}=q^{k(k-1)/2}\,(-1)^{m+1}\,\t^{m(m-1)/2}\,\t^{-km}\,\br{k-1}{k-m}_\t.\label{c}\ee
(Observe that $c_{m,k}=0$ when $m>k$.) The radius $R$ is so large that the denominators $p+q\x\x'-\x$ and $\x-1+\r\,(1-\t)$ are nonzero on and outside the contour.

We now derive an analogous formula for ASEP on $[1,\,N]$ with Bernoulli initial condition (subsequently we shall let $N\to\iy$). We do this by first taking the set $S$ in (\ref{PY}) as fixed and averaging over all $Y\subset[1,\,N]$. Then we take the sum over $S$. 

The only part of (\ref{PY}) that depends on $Y$ is $\t^{\s(S,\,Y)}$, and the probability of an initial configuration $Y$ is $\r^{|Y|}\,(1-\r)^{N-|Y|}$, so for the fixed $S\subset[1,\,N]$ we have to compute
\be\sum_{S\subset Y\subset[1,\,N]}\r ^{|Y|}\,(1-\r)^{N-|Y|}\,\t^{\s(S,\,Y)}.\label{Ysum}\ee

If $S=\{s_1,\ld,s_k\}$ then any element of $Y$ which is $> s_k$ does not affect $\s(S,Y)$, so we write $Y=Y_1\cup Y_2$ where $Y_1=\{y\in Y:y\le s_k\}$ and $Y_2=Y\bks Y_1$. Then (\ref{Ysum}) equals
\be\sum_{S\subset Y_1\subset[1,\,s_k]}\r^{|Y_1|}\,(1-\r)^{-|Y_1|}\,\t^{\s(S,\,Y_1)}\cdot \sum_{Y_2\subset(s_k,\,N]}\r^{|Y_2|}\,(1-\r)^{N-|Y_2|}.\label{Ysumprod}\ee

In the second factor the number of summands with $|Y_2|=\ell$ equals $N-s_k\choose\ell$, so the factor itself equals
\[\sum_{\ell=0}^{N-s_k} {N-s_k\choose\ell}\, \r^{\ell}\,(1-\r)^{N-\ell}=(1-\r)^{s_k}\,\sum_{\ell=0}^{N-s_k} {N-s_k\choose\ell} \r^{\ell}\,(1-\r)^{N-s_k-\ell}=(1-\r)^{s_k}.\]

For the first factor in (\ref{Ysumprod}) we set $t_i=s_i-s_{i-1}-1\ (s_0=0)$, the number of points in the gap between $s_i$ and $s_{i-1}$. Thus $s_i=t_1+\cd+t_i+i$. The exponent $\s(S,\,Y_1)$ is determined by the $j_i$ the number of points of $Y_1$ in these gaps. We have 
\[\s(S,\,Y_1)=(j_1+1)+(j_1+j_2+2)+\cd+(j_1+\cd+j_k+k),\]
the $i$th summand coming from the contribution of $s_i$ to $\s(S,\,Y_1)$. We have $|Y_1|=k+\sum j_i$ and the number of sets $Y_1$ with these $j_i$ is $\prod_i{t_i\choose j_i}$. Thus the first factor in (\ref{Ysumprod}) equals
\[\sum_{j_i\le t_i}\prod_i{t_i\choose j_i}\,\cdot\,\({\r\ov 1-\r}\)^{k+\sum j_i}\,\t^{k(k+1)/2+\sum(k-i+1)j_i}\]
\[=\t^{k(k+1)/2}\,
\({\r\ov 1-\r}\)^k\,\prod_i\(1+\t^{k-i+1}\,{\r\ov 1-\r}\)^{t_i},\]
by the multinomial theorem.
  
If we recall that $s_k=t_1+\cd+t_k+k$ we see that we have shown that (\ref{Ysum}) is equal to
\[\t^{k(k+1)/2}\,\r^k\;(1-\r)^{\,\sum t_i}\;
\prod_i\(1+\t^{k-i+1}\,{\r\ov 1-\r}\)^{t_i}\]
\be=\t^{k(k+1)/2}\,\r^k\;
\prod_i\(1-\r+\t^{k-i+1}\,\r\)^{t_i}.\label{Ysum2}\ee
This replaces the sum over $Y$ of $\t^{\s(S,\,Y)}$ in (\ref{PY}).

The rest of the integrand that depends on the elements of $S$ (not just $|S|=k$) is
\[\prod\x_i^{-s_i}=\prod\x_i^{-(t_1+\cd+t_i+i)}=\prod\x_i^{-i}\,\cdot\,\prod (\x_i\,\x_{i+1}\cd\x_k)^{-t_i}.\]
We multiply this by (\ref{Ysum2}) and then sum over all $S\subset[1,\ld,N]$ with $|S|=k$, i.e., over all $t_i\ge0$ with $\sum t_i+k\le N$. The sum over such $S$ equals
\[\t^{k(k+1)/2}\,\r^k\;\prod\x_i^{-i}\;\sum_{\sum t_i\le N-k}
\prod_i\({1-\r+\t^{k-i+1}\,r\ov \x_i\cd\x_k}\)^{t_i}.\]
Now we need only let $N\to\iy$ and we obtain the simple expression
\[\t^{k(k+1)/2}\,\r^k\;\prod_i {1\ov \x_i\cd\x_k-1+\r-\t^{k-i+1} \r}.\]

To recapitulate, we have shown that for Bernoulli initial condition on $\Z^+$,
\[\P(x_m(t)=x)=\sum_{k\ge1}\t^{k(k+1)/2}\,r^k\,c_{m,k}\]
\be\times\,\int_{|\x_i|=R}I(x,k,\x)\;
\prod_i {1\ov \x_i\cd\x_k-1+\r-\t^{k-i+1} \r}\,\prod_i d\x_i.\label{PB}\ee

\begin{center}{\bf III. A Combinatorial Identity}\end{center}

An important ingredient in \cite{TW1,TW2} was a combinatorial identity which, luckily, has a generalization which will allow us to derive the Fredholm determinant representation. The generalization is 
\[\sum_{\s\in\mathbb{S}_k}\,{\rm sgn}\,\s\,\prod_{i<j}{1\ov p+q\xs{i}\xs{j}-\xs{i}}\]
\[\times{1\ov (\xs{1}\xs{2}\cd \xs{k}-1+\r-\t^k\,\r)\,(\xs{2}\cd\xs{k}-1+\r-\t^{k-1}\,\r)\cd(\xs{k}-1+\\r-\t\,\r)}\]
\be=q^{k(k-1)/2}{\prod_{i<j}(\x_i-\x_j)\ov
\prod_i(\x_i-1+\r\,(1-\t))\,\cdot\,\prod_{i\ne j}(p+q\x_i\x_j-\x_i)}.\label{Biden}\ee

For the proof we shall try to be more general (the notation in the proof is also simpler this way) and see that we get nothing more than this identity.  We look for constants $a_1,\ld,a_N,\,b_1,\ld,b_N,\,c$ for which there is an identity

\[\sum_{\s\in\mathbb{S}_k}\,{\rm sgn}\,\s\,\prod_{i<j}{1\ov p+q\xs{i}\xs{j}-\xs{i}}\]\[\times{1\ov (\xs{1}\xs{2}\cd \xs{k}-a_k)\,(\xs{2}\cd\xs{k}-a_{k-1})\cd(\xs{k}-a_1)}\]
\[=b_1\cd b_k{\prod_{i<j}(\x_i-\x_j)\ov
\prod_i(\x_i-c)\,\cdot\,\prod_{i\ne j}(p+q\x_i\x_j-\x_i)}.\]

We try to prove the identity by induction, the truth for $k=1$ being equivalent to
\be b_1=1,\ \ a_1=c.\label{ab}\ee
Denote the left side of the proposed identity by $\ph_k(\x_1,\ld,\x_k)$ and the right side by $\ps_k(\x_1,\ld,\x_k)$.

We assume the identity is true for $k-1$. We first sum over all permutations such that $\s(1)=\l$, and then sum over $\l$. If we observe that the inequality $i<j$ becomes $j\ne i$ when $i=1$, we see that $\ph_k(\x_1,\ld,\x_k)$ equals
\[\sum_{\l=1}^k(-1)^{\l+1}{1\ov \prod_{j\ne \l}(p+q\x_\l\x_j-\x_\l)}\,\cdot\,{1\ov \x_1\x_2\cd \x_k-a_k}\,\cdot\,
\ph_{k-1}(\x_1,\ldots,\x_{\l-1},\x_{\l+1},\ldots,\x_k)\]
\[=\sum_{\l=1}^k(-1)^{\l+1}{1\ov \prod_{j\ne \l}(p+q\x_\l\x_j-\x_\l)}\,\cdot\,{1\ov \x_1\x_2\cd \x_k-a_k}\,\cdot\,
\ps_{k-1}(\x_1,\ldots,\x_{\l-1},\x_{\l+1},\ldots,\x_k),\]
by the induction hypothesis. This equals $\ps_k(\x_1,\ld,\x_k)$ if
\[\sum_{\l=1}^k(-1)^{\l+1}{1\ov \prod_{j\ne \l}(p+q\x_\l\x_j-\x_\l)}\,\cdot\,
{\ps_{k-1}(\x_1,\ldots,\x_{\l-1},\x_{\l+1},\ldots,\x_k)\ov\ps_{k}(\x_1,\cd,\x_k)}= \x_1\x_2\cd \x_k-a_k.\]

The left side is equal to
\[\sum_{\l=1}^k(-1)^{\l+1}{1\ov \prod_{j\ne \l}(p+q\x_\l\x_j-\x_\l)}\,\cdot\,{1\ov b_k}\,\cdot\,{(\x_\l-c)\,\prod_{j\ne \l}(p+q\x_\l\x_j-\x_\l)\,\cdot\,\prod_{i\ne \l}(p+q\x_i\x_\l-\x_i)\ov\prod_{\l<j}(\x_\l-\x_j)\,\cdot\,\prod_{i<\l}(\x_i-\x_\l)}\]
\be={1\ov b_k}\,\sum_{\l=1}^k{(\x_\l-c)\,\prod_{i\ne \l}(p+q\x_i\x_\l-\x_i)\ov\prod_{i\ne \l}(\x_\l-\x_i)}.\label{expr}\ee

We want (\ref{expr}) to equal $\x_1\x_2\cd \x_k-a_k$. To see what it actually is we consider
\[\int{(z-c)\,\prod_{i}(p+q\x_iz-\x_i)\ov(qz-p)\,(z-1)\,\prod_{i\ne \l}(z-\x_i)}\,dz\]
taken over a large circle. Since the integrand is $\sim q^{N-1}\prod\x_i/z$ as $z\to\iy$ the integral is equal to
\[ q^{k-1}\prod\x_i.\]
There are poles at $z=1,\ z=\t,$ and $z=\x_k\ (k=1,\ld,N)$. The residue at $z=1$ equals
\[{1-c\ov q-p}\,p^k.\]
The residue at $z=\t$ equals
\[{\t-c\ov q(\t-1)}\,q^k.\]
And the residue at $z=\x_\l$ equals
\[{(\x_k-c)\,\prod_{i\ne \l}(p+q\x_i\x_\l-\x_i)\ov\prod_{i\ne \l}(\x_\l-\x_i)}.\]
Here we used $p+q\x_\l^2-\x_\l=(q\x_\l-p)\,(\x_\l-1)$. 

These show that (\ref{expr}) equals
\[{1\ov b_k}\left[q^{k-1}\prod\x_i-{1-c\ov q-p}\,p^k-{\t-c\ov q(\t-1)}\,q^k\right]={q^{k-1}\ov b_k}\left[\prod\x_i-{1-c\ov 1-\t}\,\t^k-{\t-c\ov \t-1}\right].\]
This will equal $\x_1\x_2\cd \x_k-a_k$ when
\[b_k=q^{k-1},\ \ a_k={(1-c)\,\t^k+c-\t\ov 1-\t}.\]
These are consistent with (\ref{ab}). If we set 
$\r=(1-c)/(1-\t)$ then $a_i=1-\r+\r\,\t^i$ and we get (\ref{Biden}).

\begin{center}{\bf IV. Fredholm Determinant Representation}\end{center}

The part of the  integrand in (\ref{PY}) that is not symmetric in the $\x_i$ is
\[\prod_{i<j}{\x_j-\x_i\ov p+q\x_i\x_j-\x_i}\;\prod_i {1\ov \x_i\cd\x_k-1+\r-\t^{k-i+1} \r},\]
and using (\ref{Biden}) we see that the symmetrization of this is
\[{q^{k(k-1)/2}\ov k!}\,\prod_{i\ne j}{\x_j-\x_i\ov p+q\x_i\x_j-\x_i}\,\cdot\,\prod_i{1\ov\x_i-1+\r\,(1-\t)}.\]
Thus (\ref{PB}) may be replaced by
\[\P(x_m(t)=x)=\sum_{k\ge1}{1\ov k!}\,q^{k(k-1)/2}\,\t^{k(k+1)/2}\,c_{m,k}\]
\[\times\,\int_{|\x_i|=R} \prod_{i\ne j}{\x_j-\x_i\ov p+q\x_i\x_j-\x_i}\;(1-\prod_i\x_i)\,\cdot\,\prod_i{\r\ov \x_i-1+\r\,(1-\t)}\;\prod_i{\x_i^{x-1}\,e^{\ep(\x_i)t}\ov 1-\x_i}\,d\x_i.\]

Next we use the identity \cite{TW2} 
\[\det\left({1\ov p + q \x_i\x_j - \x_i}\right)=(-1)^k
 (pq)^{k(k-1)/2} \prod_{i\neq j}{\x_j-\x_i\ov p +q \x_i\x_j -\x_i}\,\,
\prod_i {1\ov (1-\x_i)(q\x_i-p)}\]
and (\ref{c}) to write the above as
\[\P(x_m(t)=x)=(-1)^{m+1}\t^{m(m-1)/2}\,\,\sum_{k\ge m}{(-1)^{k}\ov k!}\,\t^{(1-m)\,k}\,\br{k-1}{k-m}_\t\]
\[\times\,\int_{|\x_i|=R}\det\({q\ov p + q \x_i\x_j - \x_i}\)\,(1-\prod_i\x_i)\,\cdot\,\prod_i{\r\,(\x_i-\t)\ov \x_i-1+\r\,(1-\t)}\;\prod_i \x_i^{x-1}\,e^{\ep(\x_i)t}\,d\x_i.\]

If we sum on $x$ from $-\iy$ to $x$, which we may do if $R>1$, and recall the definition (\ref{K}) of $K(\x,\,\x')$ we see that this becomes
\[\P(x_m(t)\le x)=(-1)^{m}\,\t^{m(m-1)/2}\,\sum_{k\ge m}\t^{(1-m)\,k}\,\br{k-1}{k-m}_\t\]
\[\times\,{(-1)^{k}\ov k!}\,\int_{|\x_i|=R}\det(K(\x_i,\,\x_j))\;d\x_1\cd d\x_k.\]

The integral with its factor is the coefficient of $\la^k$ in the expansion of $\det(I-\la\,K)$ and so is equal to
\[\int \det(I-\la\,K)\,{d\la\ov\la^{k+1}}.\]
By the $\t$-binomial theorem \cite[p.26]{M} we have  for $|z|$ small enough 
\[\sum_{k\ge m} \br{k-1}{k-m}_\t z^{k}=z^m\sum_{j\ge 0} \br{m+j-1}{j}_\t z^j=z^m\,\prod_{i=0}^{m-1}{1\ov  1- \t^i z} =\prod_{j=1}^{m}{z\ov  1- \t^{m-j} z}.\]
If we set $z=\t^{1-m}\,\la\inv$ this gives for $|\la|$ large enough
\[(-1)^m\,\t^{m(m-1)/2}\,\sum_{k\ge m} \br{k-1}{k-m}_\t\,\t^{(1-m)\,k}\,\la^{-k}=\prod_{j=0}^{m-1}{1\ov 1-\la\,\t^j}.\]
Thus,
\[\P(x_m(t)\le x)=\int{\det(I-\la\,K)\ov \prod_{j=0}^{m-1}(1-\la\,\t^j)}\;{d\la\ov\la},\]
where the contour of integration encloses all the singularities of the integrand.

\begin{center}{\bf V. Asymptotics when \boldmath$p<q$}\end{center}

In this section we shall rely on \cite{TW3}. In order to quote the results of that paper we change the definition of our kernel slightly, replacing $t$ by $t/\g$, to read
\be K(\x,\,\x')=q\,{{\x}^x\,e^{\ep(\x)t/\g}\ov p+q\x\x'-\x}\,{\r\,(\x-\t)\ov\x-1+\r\,(1-\t)},\label{K1}\ee
so that now
\be\P(x_m(t/\g)\le x)=\int{\det(I-\la\,K)\ov \prod_{j=0}^{m-1}(1-\la\,\t^j)}\;{d\la\ov\la}.\label{P2}\ee
The difference between the kernel $K(\x,\,\x')$ here and \cite[(1)]{TW3} is that the present one has the extra factor
\[{\r\,(\x-\t)\ov\x-1+\r\,(1-\t)}.\]
If we go to $\e,\,\e'$ variables defined by the substitutions
\[\x={1-\t\e\ov1-\e},\ \ \ \x'={1-\t\e'\ov1-\e'},\]
we obtain the kernel
\[K_2(\e,\,\e')={\ph(\e')\ov \e'-\t\e}\]
where instead of 
\[\ph(\e)=\({1-\t \e\ov1-\e}\)^x\,e^{\left[{1\ov 1-\e}-
{1\ov 1-\t \e}\right]\,t}\]
as in \cite{TW3} we have now
\[\ph(\e)=\({1-\t \e\ov1-\e}\)^x\,e^{\left[{1\ov 1-\e}-
{1\ov 1-\t \e}\right]\,t}\,{1\ov1+\a\,\e},\]
with 
\[\a={1-\r\ov \r}.\]
(The operator now acts on a small circle about $\e=1$ described clockwise. Whether the factor is $\ph(\e')$ or $\ph(\e)$ makes no difference for the determinant.)

In \cite{TW3} we used lemmas on stability of the Fredholm determinant  to show that if
\[K_1(\e,\,\e')={\ph(\t\e)\ov \e'-\t\e},\]
then the Fredholm determinant of $K_2$ on the small circle equals the Fredholm determinant of $K_1(\e,\,\e')-K_2(\e,\,\e')$ on $\G$, a circle about zero described counterclockwise, with 1 inside and $\t\inv$ outside. We were able to do this because $\ph(\e)$ was analytic between the small circle and this contour. But now it may not be; it is only the case when $\a<1$, or equivalently $\r>1/2$. When $\r\le 1/2$ we take $\G$ to be the circle with diameter $[-\a\inv+\dl,1+\dl]$, with $\dl$ small. Then $\ph(\e)$ is analytic between the small circle and this contour, so that the argument of \cite{TW3} holds. Also holding is Proposition~4 of that paper which says that $\det(I-\la\,K_1)$ acting on $\G$ equals $\prod_{k=0}^{\iy}(1-\la\,\t^k)$. It was important for this, and will be for what comes later, that $\ph(0)=1$.

Following the discussion of \cite[\S V]{TW3} we define
\[f(\m,\,z)=\sum_{k=-\iy}^\iy{\t^{k}\ov 1-\t^{k}\m}\,z^k,\]
which is analytic for $1<|z|<\t\inv$ and extends analytically to all $z\ne0$ except for poles at the $\t^k,\ k\in\Z$,
\be\phy(\e)=\prod_{n=0}^\iy \ph(\t^n\,\e)=(1-\e)^{-x}\,e^{{\e\ov 1-\e}t}\,\prod_{n=0}^\iy{1\ov 1+\t^n\a\e},\label{phy}\ee
and
\be J(\e,\,\e')=\int{\phy(\z)\ov\phy(\e)}\,{\z^m\ov \e^{m+1}}\,{f(\m,\z/\e')\ov\z-\e}\;d\z.\label{J}\ee
(We replaced two occurrences of $\e'$ in the denominator by $\e$, which does not change the determinant.) When $\r>1/2$ all the singularities $-(\t^n\a)\inv$ lie outside the unit circle and we take the $\z$-contour to be a circle slightly larger than the unit circle; the operator $J$ then acts on a circle slightly smaller than the unit circle. If $\r\le 1/2$ we take the $\z$-contour to be the circle with diameter $[-\a\inv+\dl,1+\dl]$ and $J$ acts on the circle with diameter $[-\a\inv+2\dl,1-\dl]$, with $\dl$ small.

We obtain the formula
\be\P(x_m(t/\g)\le x)=\int \prod_{k=0}^{\iy}(1-\m\,\t^k)\,\cdot\,\det \(I+\m\, J\) \;{d\m\ov \m},\label{P3}\ee
where the integration is taken over a circle with center zero and radius larger than $\t$ but not equal to any $\t^{-k}$ with $k\ge 0$. 

Rather than go through the steps of \cite{TW3} to show that this is so, we use an analyticity argument. The right side of (\ref{P2}) extends to an analytic function of $\r$ in the neighborhood of $(0,\,1)$ (actually an entire function). So does (\ref{P3}): we need only vary the $\z$- and $\e$-contours when $|\r|\le 1/2$, making sure that all the singularities $-(\t^n\a)\inv$ lie outside the $\z$-contour and that the $\e$-contour lies slightly inside the $\z$-contour. Since the two sides agree when $\r\in(1/2,\,1]$, when the argument of \cite{TW3} goes over word-for-word,\footnote{All that was used was that $\phy$ was analytic and nonzero on the unit disc.} they agree for all $\r$.

In \cite{TW3} we did a saddle point analysis of (\ref{J}), in which
$x=c_1\,t+c_2\,t^{1/3}$, by writing the main part of the integrand as $\ps(\z)/\ps(\e)$, where
\be\ps(\z)=(1-\z)^{-c_1\,t}\,e^{{\z\ov 1-\z}t}\,\z^{\s t}.\label{ps}\ee
The saddle point is
 \[\x=-{\sqrt\s\ov 1-\sqrt\s}.\]
We deformed the $\z$-contour so it passed through $\x-t^{-1/3}$, we deformed the $\e$-contour so it passed through $\x$ (which did not change the Fredholm determinant), and then in the neighborhood of the saddle point made the substitutions
\be\e\to \x+c_3\inv\,t^{-1/3}\,\e,\ \ \ \e'\to\x+c_3\inv t^{-1/3}\,\e',\ \ \ \z\to \x+c_3\inv\,t^{-1/3}\,\z,\label{substitutions}\ee
with a certain constant $c_3$. The limiting $\e$-contour consisted of the rays from 0 to $c_3\,e^{\pm\pi i/3}$ while the limiting $\z$-contour consisted of the rays from $-c_3$ to $-c_3+c_3\,e^{\pm2\pi i/3}$.
To carry out the proof of Theorem~2 we get exactly the same statement, independently of $\r$, as long as in the deformation to the contours to pass through the saddle point we do not pass any of the sigularities from the last factor in (\ref{phy}). (The value of $\phy$ at the saddle point $\x$ drops out in the quotient and so is irrelevant.) This requires that 
\[-{\sqrt\s\ov 1-\sqrt\s}>-{1\ov\a}=-{\r\ov1-\r},\]
or $\s<\r^2$. This gives the first part of Theorem 2.

When $\s=\r^2$ we have $\x=-\a\inv$, so when $\s=\r^2+o(t^{-1/3})$ we have $\x=-\a\inv+o(t^{-1/3})$. Now to avoid passing across the singularity $-\a\inv$ we deform the contours so that the $\e$-contour passes through $\x+2\dl\,t^{-1/3}$ and the $\z$-contour passes through $\x+\dl\,t^{-1/3}$, where $\dl>0$ is arbitrary but fixed. In the neighborhood of the saddle point we make the same substitutions (\ref{substitutions}). Now the limiting $\e$-contour $\G_\e$ consists of the rays from $2\dl c_3$ to $2\dl c_3+c_3\,e^{\pm\pi i/3}$ while the limiting $\z$-contour $\G_\z$ consists of the rays from $\dl c_3$ to $\dl c_3+c_3\,e^{\pm2\pi i/3}$. The part of the infinite product from (\ref{phy}) with $n>1$ in the ratio $\phy(\z)/\phy(\e)$ is $1+o(1)$ in the limit. But the part with $n=1$,
\[{1+\a\e\ov1+\a\z}={\a+\e\ov\a+\z},\]
is different. After the substitutions it becomes
\[{c_3\,t^{1/3}\,(\a\inv+\x)+\e\ov c_3\,t^{1/3}\,(\a\inv+\x)+\z}={\e\ov\z}
+o(1),\]
since $\a\inv+\x=o(t^{-1/3})$.

The upshot is that the kernel
\[\int_{\G_\z}\int_{\G_\e}{e^{-\z^3/3+\e^3/3+y\z-x\e}\ov\z-\e}\,d\e\,d\z\]
at the end of \cite{TW3} gets replaced by
\[\int_{\G_\z}\int_{\G_\e}{e^{-\z^3/3+\e^3/3+y\z-x\e}\ov\z-\e}\,{\e\ov\z}\,d\e\,d\z.\]
This equals
\[-\left[\KA(x,\,y)+\A(x)\,\int_{-\iy}^y\A(z)\,dz\right],\]
which gives the distribution function $F_1(s)^2$. The limit is independent of, and uniform in, $\m$ and so we get the same limit when we do the integration in (\ref{P3}). This gives the second part of Theorem~2, with its final statement. 
 
When $\s>\r^2$ we replace the constant $c_1$ in (\ref{ps}) by $c_1'$ and set $x=c_1'\,t+c_2'\,s\,t^{1/2}$. Because of the choice of $c_1'$  the saddle point is at the singularity $-\a\inv$. In the neighborhood of the saddle point 
\[\ps(\z)= \ps(-\a\inv)\,e^{[-\r\inv(1-\r)^3\,(\s-\r^2)(\z+\a\inv)^2/2+O((\z+\a\inv)^3)]\,t}.\]
Since $\s>\r^2$ we see that if the $\e$-contour is the component of the level curve $|\ps(\e)/\ps(-\a\inv)|=1-\dl$ with the saddle point outside and a small indentation to the left of $\e=1$, and the $\z$-contour is expanded to the component of the level curve $|\ps(\z)/\ps(-\a\inv)|=1-2\dl$ with the saddle point inside and a small indentation to the right of $\z=1$ then the norm of the operator represented by the new $\z$-integral is exponentially small. When we expanded the $\z$-contour we picked up the residue
\[{\ps(-\a\inv)\ov\e\,\ps(\e)}\,f(\m,-(\a\e')\inv)\,\({1-\e\ov1+\a\inv}\)^{c_2'\,t^{1/2}}\,\prod_{n=1}^\iy{1+\t^n\a\e\ov1-\t^n}.\]
from the singularity. Therefore with exponentially small error $\det(I+\m\,J)$ equals
\[1+\m\int{\ps(-\a\inv)\ov\ps(\e)}\,\e\inv\,f(\m,-(\a\e)\inv))\,\({1-\e\ov1+\a\inv}\)^{c_2'\,s\,t^{1/2}}\,\prod_{n=1}^\iy{1+\t^n\a\e\ov1-\t^n}\,d\e,\]
assuming this is bounded away from zero  as $t\to\iy$.

We evaluate this integral by steepest descent. The steepest descent curve passes through (i.e., just to the right of) the saddle point $-\a\inv$ vertically, loops around the origin and has an inward-facing cusp at $\e=1$. In deforming the previous $\e$-curve to this one we pass no singularities of the integrand. In the neighborhood of the saddle point we have 
\[\e\inv\,f(\m,-(\a\e)\inv))={\m\inv\ov \e+\a\inv}+O(1),\]
\[\({1-\e\ov1+\a\inv}\)^{c_2'\,t^{1/2}}=e^{-[(1-\r)\,c_2'\,s\,(\e+\a\inv)+O((\e+\a\inv)^2)]\,t^{1/2}},\]
and the last factor in the integrand is $1+O(\e+\a\inv)$.
Therefore, if we deform to the steepest descent contour and make the replacement $\e\to \e-\a\inv$ we see that we get the limit
\[1-\int_{-i\iy+0}^{i\iy+0}e^{\,\r^{-2}(1-\r)^3\,(\s-\r^2)\e^2/2-(1-\r)\,c_2'\,s\,\e}\,{d\e\ov\e}=1-\int_{-i\iy+0}^{i\iy+0}e^{\e^2/2-\,s\,\e}\,{d\e\ov\e}=G(s),\]
where we used $\r^{-2}(1-\r)^3\,(\s-\r^2)=(1-\r)^2\,{c_2'}^2$. (Of course this is how $c_2'$ was chosen.) 
Again, this limit being independent of $\m$, the integral (\ref{P3}) has the same limit. This completes the proof of Theorem~2.
\sp

As was seen in \cite{TW4} in the case $\r=1$, Theorem~3 becomes a corollary of Theorem~2 if we use the fact
\[\P(\T(x,\,t)\ge m)=\P(x_m(t)\le x).\]
Here is the argument for parts one and two. If we set $x=vt$ and let $m=a_1\,t+s\,a_2\,t^{1/3}$, then the left side becomes
\[\P\(\T(vt,\,t)\ge a_1\,t+s\,a_2\,t^{1/3}\).\]
We will be able to apply Theorem~2 with $s$ replaced by some $s'$ if when
\be\s=a_1+s\,a_2\,t^{-2/3}\label{siga}\ee
then $a_1$ and $a_2$ are such that, in the notation of Theorem~2,
\[c_1 \,t+c_2 \,s'\,t^{1/3}=vt.\]
(And $\s$ is in the range for which the conclusion is applicable.) Substituting what $c_1 $ and $c_2 $ are and dividing by $t$ give
\[-1+2\sqrt{\s}+\s^{-1/6}\,(1-\sqrt{\s})^{2/3}\,s'\,t^{-2/3}=v,\]
and solving for $\s$ gives
\[\s=(1+v)^2/4-2^{-4/3}(1-v^2)^{2/3}\,s'\,t^{-2/3}+O(t^{-4/3}).\]
This becomes (\ref{siga}), with $a_1$ and $a_2$ as in the statement of Theorem~3, if $s'=-s+O(t^{-2/3})$. 

Applying the first two parts Theorem~2, with their uniformity statements, gives the first two parts of the theorem. (Observe that $\s=\r^2+O(t^{-2/3})$ when $v=2\r-1$, which is why we could not simply let $\s=\r^2$ in part two of Theorem~2.) The third part follows in the same way from the third part of Theorem~2. 

\begin{center}{\bf Acknowledgments}\end{center}

The authors thank Herbert Spohn and Patrik Ferrari for useful references. Special thanks go to Ivan Corwin for a fruitful discussion related to paper \cite{BC}.

This work was supported by the National Science Foundation through grants DMS-0553379 and DMS-0906387 (first author) and DMS-0552388 and DMS-0854934 (second author).

\end{document}